\documentclass{amsart}
\usepackage{amsmath,amssymb,latexsym}
\usepackage{diagrams}
\usepackage[mathscr]{eucal}

\newcommand\Pic[1]{\hbox{\rm Pic(}#1\hbox{\rm )}}

\newcommand\sO{{\mathcal O}}

\newcommand\bP{{\mathbb P}}

\newcommand\bC{{\mathbb C}}

\def\rank{\operatorname{rank}}

\def\lim{\operatorname{lim}}

\def\Leff{\operatorname{Leff}}

\def\dlog{\operatorname{dlog}}

\def\NS{\operatorname{NS}}
\def\Num{\operatorname{Num}}

\def\Pic{\operatorname{Pic}}

\def\id{\operatorname{id}}

\newtheorem{theorem}{Theorem}[section]
\newtheorem{lemma}[theorem]{Lemma}
\newtheorem{corollary}[theorem]{Corollary}

\theoremstyle{definition}
\newtheorem{definition}[theorem]{Definition}

\theoremstyle{remark}

\newtheorem{remarks}[theorem]{Remarks}

\numberwithin{equation}{section}

\begin{document}

\title{On the normal bundle of submanifolds of $\mathbb P^n$}
\author{Lucian B\u adescu}
\address{Dipartimento di Matematica, Universit\`a degli Studi di Genova,
Via Dodecaneso 35, 16146 Genova, Italy}
\email{badescu@dima.unige.it}

\date{June 30, 2006}



\begin{abstract} Let $X$ be a submanifold of dimension $d\geq 2$ of the complex projective space $\mathbb P^n$. We prove results of the following type.   i) If $X$ is irregular and $n=2d$ then the normal bundle $N_{X|\mathbb P^n}$ is indecomposable. ii)  If $X$ is irregular, $d\geq 3$ and $n=2d+1$ then $N_{X|\mathbb P^n}$ is not the direct sum of two vector bundles of rank $\geq 2$. iii) If $d\geq 3$, $n=2d-1$ and $N_{X|\mathbb P^n}$ is decomposable then the natural restriction map $\Pic(\mathbb P^n)\to\Pic(X)$ is an isomorphism (and in particular, if $X=\mathbb P^{d-1}\times\mathbb P^1$ embedded Segre in $\mathbb P^{2d-1}$ then $N_{X|\mathbb P^{2d-1}}$ is indecomposable). iv) Let $n\leq 2d$ and $d\geq 3$, and assume that $N_{X|\mathbb P^n}$ is a direct sum of line bundles; if $n=2d$ assume furthermore that $X$ is simply connected and $\mathscr O_X(1)$ is not divisible in $\Pic(X)$. Then $X$ is a complete intersection. These results follow from Theorem \ref{exact5} below together with Le Potier vanishing theorem. The last statement also uses a criterion of Faltings for complete intersection. In the case when $n<2d$ this fact was proved by M. Schneider in 1990 in a completely different way. \end{abstract}

\maketitle

\section*{Introduction} It is well known that if $X$ is a submanifold of the complex projective space $\mathbb P^{n}$ 
($n\geq 3$) of dimension $d>\frac{n}{2}$ then a topological result of of Lefschetz type, due to Barth and Larsen (see  
\cite{[L]}, \cite{[Ba]}), asserts that the canonical restriction maps
$$H^i(\mathbb P^{n},\mathbb Z)\to H^i(X,\mathbb Z)$$ 
are  isomorphisms for $i\leq 2d-n$, and injective with torsion-free cokernel, for $i=2d-n+1$. As a consequence, the restriction map $$\Pic(\mathbb P^{n})\to\Pic(X)$$ is an isomorphism if $d\geq\frac{n+2}{2}$, and injective with torsion-free cokernel if $n=2d-1$. In particular, if $d>\frac{n}{2}$ then the class of $\mathscr O_X(1)$ is not divisible in  $\Pic(X)$. 

In this paper, in the spirit of Barth-Larsen theorem, we are going to say something about the normal bundle $N_{X|\mathbb P^{n}}$ of submanifolds $X$ of dimension $d\geq 3$ of $\mathbb P^{n}$. Specifically, we shall prove that if $X$ is a submanifold of dimension $d\geq 3$ of $\mathbb P^{2d-1}$ whose the normal bundle $N_{X|\mathbb P^{2d-1}}$ is decomposable then the restriction map $\Pic(\mathbb P^{2d-1})\to\Pic(X)$ is an isomorphism (see Theorem \ref{2n+1}, (1) below). In particular, the normal bundle of the image of the Segre embedding $\mathbb P^{d-1}\times\mathbb P^1\hookrightarrow\mathbb P^{2d-1}$ is indecomposable for every $d\geq 3$. This result suggests that the decomposability of the normal bundle of a given submanifold $X$ of $\mathbb P^n$ of dimension $d\geq 3$ should yield strong geometrical constraints on $X$. For illustration, see Theorem \ref{2n+1} and its corollaries. For example Theorem \ref{2n+1}, (3) asserts that every submanifold of $\mathbb P^n$ of dimension $d\geq 3$, whose normal bundle is a direct sum of line bundles, is regular and has $\Num(X)$ is isomorphic to $\mathbb Z$ (here $\Num(X):=\Pic(X)/{\text{numerical equivalence}}$); moreover, if either $2d>n$, or if $n=2d$, $X$ is simply connected and $\mathscr O_X(1)$ is not divisible in $\Pic(X)$, then $X$ is a complete intersection (if $d>\frac{n}{2}$ this result was first proved, in a different way, by M. Schneider in \cite{[S]}). Another result (Theorem \ref{ci11}) asserts the following: (1) the normal bundle of any irregular submanifold of dimension $d\geq 2$ in $\mathbb P^{2d}$ is indecomposable; (2) the normal bundle of any irregular submanifold of dimension $d\geq 3$ in $\bP^{2d+1}$ is not the direct sum of two vector bundles of rank $\geq 2$.

Although these kind of results seem to be completely new, the idea behing their proofs is surpringly simple. Our basic technical result (Theorem \ref{exact5}) asserts that for every submanifold $X\subset\mathbb P^n$ of dimension $d\geq 2$ the irregularity of $X$ is equal to $h^1(N_{X|\mathbb P^n}^{\vee})$, and the rank of the N\'eron-Severi group of $X$ is $\leq 1+h^2(N_{X|\mathbb P^n}^{\vee})$. This theorem and a systematic use of Le Potier vanishing theorem yield the proofs of most of the results of this paper. Certain applications of Theorem 
\ref{exact5} will also make use of a criterion of Faltings \cite {[F2]} for complete intersection.

The general philosophy according to which there is a close relationship between topological Barth-Lefschetz theorems (see \cite{[Ba]}, \cite{[Ha1]}) and vanishing results involving the conormal bundle of the variety in question is not new. For instance, Faltings showed in \cite{[F3]} that for any $d$-fold $X$ in 
$\mathbb P^n$ the following implication holds:
$$H^q(\mathbb P^n,X;\mathbb C)=0\;\;\text{for $q\leq 
2d-n+1$}\;\;\;\Longrightarrow $$
$$\Longrightarrow\;H^q(S^k(N_{X|\mathbb P^n}^{\vee}))=0\;\;\text{for $q\leq 2d-n$ and $k\geq 1$}.$$

Conversely, Schneider and Zintl proved the following vanishing result (see  \cite{[SZ]}) 
\begin{equation}\label{sz}
H^q(S^k(N_{X|\mathbb P^n}^{\vee})(-i))=0\;\;\text{for $q\leq 2d-n$, $k\geq 1$ and $i\geq 0$},
\end{equation}
without using Barth-Lefschetz theorem (here $E^{\vee}$ denotes the dual of a vector bundle $E$). Moreover they showed that \eqref{sz} implies Barth-Lefschetz theorem, i.e. $H^q(\mathbb P^n,X;\mathbb C)=0$, $\forall q\leq 2d-n+1$. Finally, we mention the papers \cite{[AO]}, \cite{[E]} and \cite{[PPS]} to illustrate how certain vanishings of the cohomology involving the normal bundle may have interesting geometric consequences concerning small codimensional submanifolds of $\mathbb P^n$.

\section{Some known results and background material}

All varieties considered here are defined over the field $\mathbb C$ of complex numbers. By a submanifold of $\mathbb P^n$ we mean a smooth closed irreducible subvariety of $\mathbb P^n$. The rest of the terminology and notation used throughout this paper are standard. In particular, for every projective variety $X$ one defines:

-- $\Pic^0(X)$ (resp. $\Pic^{\tau}(X)$) as the subgroup of $\Pic(X)$ of all isomorphism classes of line bundles on $X$ which are algebraically (resp. numerically) equivalent to zero. One has $\Pic^0(X)\subseteq\Pic^{\tau}(X)$ and a result of Matsusaka asserts that $\Pic^{\tau}(X)/\Pic^{0}(X)$ is a finite group (see e.g. \cite{[K]}).

--  $\NS(X):=\Pic(X)/\Pic^{0}(X)$ (the N\'eron-Severi group of $X$) 
and $\Num(X):=\Pic(X)/\Pic^{\tau}(X)$.

\medskip

The main tool used in this paper is the following generalization of Kodaira vanishing theorem due to Le Potier:

\begin{theorem}[Le Potier vanishing theorem \cite{[LP]}]\label{(1.8.5)} Let $E$ be an ample vector bundle of rank $r$ on a complex projective manifold $X$ of dimension $d\geq 2$. Then $H^i(E^{\vee})=0$ for every $i\leq d-r$.\end{theorem}

We shall also need the following criterion of Faltings for complete intersection:

\begin{theorem}[Faltings \cite{[F2]}]\label{ci1} Let $X$ be a submanifold of $\mathbb P^n$ such that there is a surjection $\bigoplus\limits_{i=1}^p\mathscr O_X(a_i)\twoheadrightarrow N_{X|\mathbb P^n}$ for some positive integers $a_i$, $i=1,...,p$. If $p\leq\frac{n}{2}$ then $X$ is a complete intersection.\end{theorem}

\medskip
\medskip

Now we shall need a definition and some preliminary general results that shall be needed in the sequel. Let
\begin{equation}\label{A0}
0\to E_1\to E \to E_2\to 0
\end{equation}
be an exact sequence of vector bundles on a projective manifold $X$. If $E$ is ample, it is well known that $E_2$ is also ample, but this is not longer true in general for $E_1$. 

\begin{definition}\label{A1} Let $E$ be an ample vector bundle of rank $r\geq 2$ on a projective manifold $X$. Let
$p$ a natural number such that $1\leq p\leq\frac{r}{2}$. We say that $E$ satisfies condition $\mathscr A_p$ if there exists no exact sequence of the form \eqref{A0} with $E_1$ and $E_2$ ample vector bundles on $X$ of rank $\geq p$.
\end{definition}

Clearly, $\mathscr A_1\Rightarrow\mathscr A_2\Rightarrow\cdots$. On the other hand, if an ample vector bundle $E$ satisfies condition $\mathscr A_1$ then $E$ is indecomposable, i.e. $E$ cannot be written as $E=E_1\oplus E_2$, with $E_1$ and 
$E_2$ vector bundles of rank $\geq 1$. We are going to apply Definition \ref{A1} to the normal bundle $N_{X|\bP^n}$ of a submanifold $X$ of $\mathbb P^n$.

\medskip

First we note the following general essentially well known fact (see \cite{[G]} and \cite{[E]} for some special cases):

\begin{lemma}\label{r2} Assume that $X$ is a submanifold of dimension $d\geq 1$ of the projective space $\mathbb P^n$, such that the projection $\pi_P\colon\mathbb P^n\setminus\{P\}\to\mathbb P^{n-1}$ of center a general point $P\not\in X$ defines a biregular isomorphism $X\cong X':=\pi_P(X)$. Then there exists a canonical exact sequence
\begin{equation}\label{e3}
0\to\mathscr O_X(1)\to N_{X|\mathbb P^n}\to N_{X'|\mathbb P^{n-1}}\to 0. \end{equation}
In particular, under the above hypotheses, the normal bundle $N_{X|\mathbb P^n}$ does not satisfy condition $\mathscr A_1$.
\end{lemma}

The proof is standard and we omit it.
We also notice the following well known and simple fact:

\begin{lemma}\label{A3} Let $X$ be a submanifold of $\mathbb P^n$ of dimension $d\geq 1$, with $n\geq 2d+1$. Then
there exists an exact sequence of vector bundles on $X$ of the form
$$0\to\mathscr O_X(1)\to N_{X|\mathbb P^n}\to E\to 0.$$
In particular, if $n\geq 2d+1$, then $N_{X|\mathbb P^n}$ does not satisfy condition $\mathscr A_1$.\end{lemma}

\proof From the Euler sequence restricted to $X$ we get a surjection $\mathscr O_X^{\oplus n+1}\twoheadrightarrow  N_{X|\mathbb P^n}(-1)$, i.e. $N_{X|\mathbb P^n}(-1)$ is generated by its global sections. The hypothesis that $n\geq 2d+1$
translates into $\rank(N_{X|\mathbb P^n}(-1))\geq d+1$. Then by a well known theorem of Serre, there exists a nowhere vanishing section $s\in H^0(N_{X|\mathbb P^n}(-1))$. Thus $s$ yields an exact sequence 
$$0\to\mathscr O_X\to N_{X|\mathbb P^n}(-1)\to F\to 0$$
of vector bundles. Twisting by $\mathscr O_X(1)$ we get the desired exact sequence.\qed

\section{A general result on submanifolds in $\mathbb P^n$}

The aim of this section is to prove the following:

\begin{theorem}\label{exact5} Let $X$ be a projective submanifold  of dimension $d\geq 2$ of $\mathbb P^n$. Then: 
\begin{enumerate}
\item $h^1(\mathscr O_X)=h^1(N_{X|\mathbb P^n}^{\vee})$. In particular, $X$ is regular if and only if 
$H^1(N_{X|\mathbb P^n}^{\vee})=0$.
\item For every $i$ such that $2\leq i\leq d$ one has
$h^{i-1}(\Omega_X^1)\leq h^{i-2}(\mathscr O_X)+h^i(N_{X|\mathbb P^n}^{\vee})$. In particular, $h^1(\Omega^1_X)\leq 1+h^2(N_{X|\mathbb P^n}^{\vee})$. If $d\geq 3$ and $H^1(\mathscr O_X)=0$ the latter inequality becomes equality.  
\item $\rank\,\Num(X)\leq 1+h^2(N_{X|\mathbb P^n}^{\vee})$. In particular, if 
$H^2(N_{X|\mathbb P^n}^{\vee})=0$ then $\Num(X)\cong\mathbb Z$.\end{enumerate}
\end{theorem}

\proof Much of the geometric information about the embedding $X\subseteq\mathbb P^n$ is contained in the following commutative diagram with exact rows and columns:
\begin{diagram}
&      &0&             &0\\
&     &\dTo&            &\dTo\\
&      &\mathscr O_X&\rTo{\id} &\mathscr O_X\\
&      &\dTo&             &\dTo \\
0&\rTo &F&\rTo        &\mathscr O_X(1)^{\oplus n+1}&\rTo  &N_{X|\mathbb 
P^{n}}&\rTo &0\\
&      &\dTo&              &\dTo&                  &\dTo_{\id}\\
0&\rTo &T_X&\rTo         &T_{\mathbb P^{n}}|X&\rTo & N_{X|\mathbb 
P^n}&\rTo &0\\
&      &\dTo&             &\dTo\\
&      &0&              &0\\
\end{diagram} 
in which the last row is the normal sequence of $X$ in $\mathbb P^n$ and the middle column is the Euler sequence of $\mathbb P^n$ restricted to $X$. Analogous diagrams have been already used in literature in a crucial way to prove some results of projective geometry (see e.g. \cite{[E]}, or \cite{[B]}, pages 7 and 25). Notice that the 
sheaf $F^{\vee}$ coincides to ${\mathscr P}^1(\mathscr O_X(1))(-1)$, where ${\mathscr P}^1(\mathscr O_X(1))$ is the sheaf of first order principal parts of $\mathscr O_X(1)$. 

Dualizing the middle row and the first column we get the exact sequences
\begin{equation}\label{exact1}
0\to N_{X|\mathbb P^n}^{\vee}\to\mathscr O_X(-1)^{\oplus n+1}\to F^{\vee}\to 0,\end{equation}
\begin{equation}\label{exact2}
0\to\Omega^1_X\to F^{\vee}\to\mathscr O_X\to 0.\end{equation}
Then the cohomology of \eqref{exact1} yields the exact sequence
$$H^{i-1}(\mathscr O_X(-1)^{\oplus n+1})\to  H^{i-1}(F^{\vee})\to H^i(N_{X|\mathbb P^n}^{\vee})\to H^i(\mathscr O_X(-1)^{\oplus n+1}).$$
By Kodaira vanishing theorem the first (resp. the last) space is zero for $1\leq i\leq d$ (resp. for $1\leq i\leq d-1$). Thus 
\begin{equation}\label{exact3}
h^{i-1}(F^{\vee})\leq h^{i}(N_{X|\mathbb P^n}^{\vee}),\;\;\text{for all $1\leq i\leq d$, with equality for $1\leq i\leq d-1$}.
\end{equation}

On the other hand, the exact sequence \eqref{exact2} does not split. Indeed, the dual Euler sequence 
$$0\to\Omega^1_{\bP^n}\to\sO_{\bP^n}(-1)^{\oplus n+1}\to\sO_{\bP^n}\to 0$$
corresponds to a generator of the one-dimensional $\bC$-vector space $H^1(\bP^n,\Omega^1_{\bP^n})$.
Then the exact sequence \eqref{exact2} corresponds to the image of this generator under the composite map
$$H^1(\bP^n,\Omega^1_{\bP^n})\to H^1(X,\Omega^1_{\bP^n}|X)\to H^1(X,\Omega^1_X), $$
which is known to be non zero (otherwise the class of  $\sO_X(1)$ would be zero in $H^1(X,\Omega^1_X)$).

Then the cohomology of \eqref{exact2} yields the exact sequence
$$0\to H^0(\Omega^1_X)\to H^0(F^{\vee})\to H^0(\mathscr O_X)\to H^1(\Omega^1_X)\to H^1(F^{\vee})\to H^1(\mathscr O_X).$$
Since the the exact sequence \eqref{exact2} does not split and $H^0(\mathscr O_X)=\mathbb C$, we get the following isomorphism and exact sequence:
\begin{equation}\label{exact4}
H^0(\Omega^1_X)\cong H^0(F^{\vee})\;\;\text{and}\;\;0\to  H^0(\mathscr
O_X)\to H^1(\Omega^1_X)\to H^1(F^{\vee})\to H^1(\mathscr O_X),\end{equation}
Moreover, for every $3\leq i\leq d$ we have the cohomology sequence
\begin{equation}\label{exact12}
\cdots\to H^{i-2}(\mathscr O_X)\to H^{i-1}(\Omega^1_X)\to H^{i-1}(F^{\vee})\to\cdots.\end{equation}

Now we prove $(1)$. From \eqref{exact3} we get $h^1(N_{X|\mathbb P^n}^{\vee})= h^0(F^{\vee})$, and using the isomorphism of
\eqref{exact4}, this equality becomes $h^1(N_{X|\mathbb P^n}^{\vee})=h^0(\Omega_X^1)$. Then one concludes by Serre's GAGA and the Hodge symmetry (which yield $h^0(\Omega_X^1)=h^1(\mathscr O_X)$). 

\medskip

$(2)$ From the exact sequence \eqref{exact12} we get $h^{i-1}(\Omega^1_X)\leq h^{i-2}(\mathscr O_X)+h^{i-1}(F^{\vee})$, and from \eqref{exact3}, $h^{i-1}(F^{\vee})\leq h^i(N_{X|\mathbb P^n}^{\vee})$, whence we get the first part. The second part follows from the first one and from the exact sequence of \eqref{exact4}.

\medskip

$(3)$ follows from the last part (2) and from the following standard argument (cf. \cite{[Ha]}, \cite{[EGPS]} and 
\cite{[Br]}).
Consider the (logarithmic derivative) map 
$$\dlog:\Pic(X)\to H^1(\Omega^1_X)$$ 
defined in the following way. If $Z$ is a scheme let us denote by $\mathscr O^*_Z$ the sheaf of multiplicative groups of all nowhere vanishing functions on $Z$. If $[L]\in\Pic(X)$ is represented by the $1$-cocycle $\{\xi_{ij}\}_{i,j}$ of $\mathscr O_X^*$ with respect to an affine covering $\{U_i\}_i$ of $X$ (with $\xi_{ij}\in\Gamma(U_i\cap U_j,\mathscr O_X^*)$), then
$\dlog(\{\xi_{ij}\})$ is by definition the cohomology class of the $1$-cocycle $\{\frac{d\xi_{ij}}{\xi_{ij}}\}_{i,j}$ of $\Omega^1_X$. Since $\dlog(\Pic^0(X))=0$ the map $\dlog$ yields the map $\dlog:\NS(X)=\Pic(X)/\Pic^0(X)\to H^1(\Omega^1_X)$. Moreover, by a result of Matsusaka, $\Pic^{\tau}(X)/\Pic^0(X)$ is a finite subgroup of $\NS(X)$ (see e.g. \cite{[K]}). Since the underlying abelian group of the $\mathbb C$-vector space $H^1(\Omega^1_X)$ is torsion-free it 
follows that $\dlog(\Pic^{\tau}(X))=0$. In other words, there is a unique map $\alpha:\Num(X)\to H^1(\Omega^1_X)$ such that $\dlog=\alpha\circ\beta$, where $\beta:\Pic(X)\to\Num(X)$ is the canonical surjection. Then (3) follows from the following general well known fact:

\medskip

{\em Claim:} $\alpha$ induces an injective map $\alpha_{\mathbb C}:\Num(X)\otimes_{\mathbb Z}\mathbb C\to H^1(\Omega^1_X)$.

\medskip

This proves the theorem.\qed

\begin{remarks}\label{obstruction} i) Let $(X,N)$ be a pair consisting of a projective manifold $X$ of dimension $d\geq 2$ and a vector bundle $N$ of rank $r$ on $X$. We may ask the following question: under which conditions there exists a projective embedding of $X\hookrightarrow\mathbb P^{d+r}$ such that $N_{X|\mathbb P^{d+r}}\cong N$? Theorem \ref{exact5} provides necessary conditions for $(X,N)$ for the existence of such a projective embedding. Specifically, the irregularity $h^1(\mathscr O_X)$ should be equal to $h^1(N^{\vee})$ and the inequalities $h^{i-1}(\Omega_X^1)\leq h^{i-2}(\mathscr O_X) +h^i(N^{\vee})$ should hold for every $i$ such that $2\leq i\leq d$.\end{remarks}

ii) A simple consequence of Theorem \ref{exact5} is the following weak form of Barth-Larsen theorem: if $d\geq\frac{n+2}{2}$ then $\Pic(X)\cong\mathbb Z$.
Indeed, the ineguality $d\geq\frac{n+2}{2}$ and Le Potier vanishing theorem imply that $H^i(N_{X|\mathbb P^n}^{\vee})=0$ for $i=1,2$. Moreover, Fulton-Hansen connectedness theorem (see \cite{[FH]}, or also 
\cite{[B]}, Theorem 7.4) implies that $X$ is also algebraically simply connected and that $\Pic^{\tau}(X)=0$. Thus by Theorem \ref{exact5}, $\Num(X)=\Pic(X)\cong\mathbb Z$.

\section{Applications}

\medskip

The first application of Theorem \ref{exact5} regards the normal bundle of some irregular $d$-folds in $\bP^n$. 
By Barth-Lefschetz theorem, every  $d$-fold of $\bP^n$, with $d>\frac{n}{2}$, is regular. So the first cases
to consider are $n=2d$ and $n=2d+1$. Precisely, we have the following:

\begin{theorem}\label{ci11} Let $X$ be a submanifold of dimension $d\geq 2$ of $\mathbb P^{n}$. Then:
\begin{enumerate}
\item Assume that $X$ is irregular and $n=2d$, e.g. an elliptic scroll of dimension $d\geq 2$ in $\mathbb P^{2d}$ $($by \cite{[I]} such scrolls do exist for every $d\geq 2)$. Then $N_{X|\mathbb P^{2d}}$ satisfies condition 
$\mathscr A_1$ of Definition $\ref{A1}$, and in particular, $N_{X|\mathbb P^{2d}}$ is indecomposable.
\item Assume that $d\geq 3$ and $n=2d+1$. Then $N_{X|\mathbb P^{2d+1}}$ satisfies condition $\mathscr A_2$, but never satisfies $\mathscr A_1$. In particular, $N_{X|\mathbb P^{2d+1}}$ cannot be the direct sum of two vector bundles of rank $\geq 2$.\end{enumerate}\end{theorem}

\proof (1) Assume that there is an exact sequence of the form 
$$0\to E_1\to N_{X|\mathbb P^{2d}}\to E_2\to 0,$$
with $E_1$ and $E_2$ ample vector bundles on $X$ of ranks $\geq 1$. Dualizing and taking cohomology we get
$$H^1(E_2^{\vee})\to H^1(N_{X|\mathbb P^{2d}}^{\vee})\to H^1(E_1^{\vee}).$$
Since $E_1$ and $E_2$ are both ample of rank $\leq d-1$ on the projective $d$-fold $X$, the first and the third space are zero by Le Potier vanishing theorem. It follows that $H^1(N_{X|\mathbb P^{2d}}^{\vee})=0$. Then by Theorem \ref{exact5}, (1), $H^1(\mathscr O_X)=0$. But this is impossible because $X$ was an irregular manifold by hypothesis.

\medskip

(2) We proceed similarly as in case (1). First, the fact that $N_{X|\mathbb P^{2d+1}}$ does not satisfy condition 
$\mathscr A_1$ follows from Lemma \ref{A3}. To check condition $\mathscr A_2$, assume that there exists an exact sequence of the form 
$$0\to E_1\to N_{X|\mathbb P^{2d+1}}\to E_2\to 0,$$ 
with $E_1$ and $E_2$ ample vector bundles on $X$ of rank $\geq 2$; in particular, $E_1$ and $E_2$ have both rank $\leq d-1$. Thus by Le Potier vanishing theorem  $H^1(E_1^{\vee})=H^1(E_2^{\vee})=0$, whence the cohomology sequence
$$H^1(E_2^{\vee})\to H^1(N_{X|\mathbb P^{2d+1}}^{\vee})\to H^1(E_1^{\vee})$$
yields $H^1(N_{X|\mathbb P^{2d+1}}^{\vee})=0$. Then Theorem \ref{exact5}, (1), implies $H^1(\mathscr O_X)=0$, a contradiction. \qed

\medskip

As a second application of Theorem \ref{exact5} we have the following:

\begin{theorem}\label{2n+1} Let $X$ be a submanifold of dimension $d\geq 2$ of $\mathbb P^{n}$. Then:
\begin{enumerate} 
\item Assume $d\geq 3$ and $n=2d-1$. If $N_{X|\mathbb P^{2d-1}}$ does not satisfy condition 
$\mathscr A_1$ $($e.g.
if $N_{X|\mathbb P^{n}}$ is decomposable$)$ then $\Pic(X)\cong\mathbb Z[\mathscr O_X(1)]$.
\item Assume $d\geq 4$ and $n=2d$. If $N_{X|\mathbb P^{2d}}$ does not satisfy condition $\mathscr A_2$ 
$($e.g. if $N_{X|\mathbb P^{n}}$ is the direct sum of two vector bundles of rank $\geq 2)$ then $\Num(X)\cong\mathbb Z$.
\item Assume that $N_{X|\mathbb P^{n}}$ is direct sum of line bundles.  Then $H^1(\mathscr O_X)=0$, and if $d\geq 3$,  
$\Num(X)\cong\mathbb Z$.\end{enumerate}\end{theorem}

\proof (1) Assume that there is an exact sequence of the form 
\begin{equation}\label{A4}
0\to E_1\to N_{X|\mathbb P^{2d-1}}\to E_2\to 0,\end{equation} with $E_1$ and $E_2$ ample vector bundles on $X$ of rank $\geq 1$ (in particular, $E_1$ and $E_2$ are both of rank $\leq d-2$). Then the cohomology sequence of the dual of 
\eqref{A4}
$$H^1(E_2^{\vee})\to H^1(N_{X|\mathbb P^{2d-1}}^{\vee})\to H^1(E_1^{\vee})$$
and  yield $H^2(N_{X|\mathbb P^{2d-1}}^{\vee})=0$. Then by Theorem \ref{exact5}, (3), $\Num(X)\cong\mathbb Z$. Now, Fulton-Hansen connectedness theorem (see 
\cite{[FH]}) implies that $\Pic^{\tau}(X)=0$, i.e. $\Num(X)=\Pic(X)$, whence $\Pic(X)\cong\mathbb Z$. On the other hand, the results of Barth-Larsen \cite{[L]} or of Faltings \cite{[F1]} (see also \cite{[B]}, Theorem 10.3 and Proposition 
10.10) imply that $\mathscr O_X(1)$ is not divisible in $\Pic(X)$, whence $\Pic(X)\cong\mathbb Z[\mathscr O_X(1)]$.

\medskip

The proof of part (2) is completely similar. In fact, assume that there exists an exact sequence
$$0\to E_1\to N_{X|\mathbb P^{2d}}\to E_2\to 0,$$
with $E_1$ and $E_2$ ample vector bundles of rank $\geq 2$. Since $\rank(E_1)+\rank(E_2)=d$, it follows that $\rank(E_1),\rank(E_2)\leq d-2$. Therefore, by Le Potier vanishing theorem,
$H^2(E_1^{\vee})=H^2(E_2^{\vee})=0$. Thus the cohomology sequence
$$ H^2(E_2^{\vee})\to H^2(N_{X|\mathbb P^{2d}}^{\vee})\to H^2(E_1^{\vee})$$
yields $H^2(N_{X|\mathbb P^{2d}}^{\vee})=0$. Then the conclusion follows from Theorem 
\ref{exact5}, (3).

\medskip

(2) The hypotheses and the Kodaira vanishing theorem imply $H^1(N_{X|\mathbb P^n}^{\vee})=0$ if $d\geq 2$ and also $H^2(N_{X|\mathbb P^n}^{\vee})=0$ if $d\geq 3$. Thus by  Theorem \ref{exact5} we get $H^1(\mathscr
O_X)=0$ if $d\geq 2$ and $\Num(X)\cong\mathbb Z$ if $d\geq 3$.\qed

\medskip

Here are some corollaries of Theorem \ref{2n+1}:

\begin{corollary}\label{segre1} The normal bundle $N$ of the Segre embedding $i\colon\mathbb P^{d-1}\times\mathbb P^1\hookrightarrow\mathbb P^{2d-1}$ $(d\geq 3)$ satisfies condition $\mathscr A_1$. In particular, 
$N$ is indecomposable.\end{corollary}

\proof Direct consequence of Theorem \ref{2n+1}, (1).\qed

\begin{corollary}\label{ci5} The normal bundle of $X=\mathbb P^1\times\mathbb P^1\times\mathbb P^1$ in $\mathbb P^7$ 
$($via the Segre embedding$)$ is not a direct sum of line bundles.\end{corollary}

\proof Direct consequence of Theorem \ref{2n+1}, (3).\qed

\begin{corollary}\label{segre2} Let $N$ be the normal bundle of the Segre embedding $i\colon\mathbb P^2\times\mathbb P^2\hookrightarrow\mathbb P^8$. Then $N$ satisfies condition $\mathscr A_2$.
Furthermore, there exists an exact sequence of the form
$$0\to \mathscr O_{\mathbb P^2\times\mathbb P^2}(1,1)\to N\to N'\to 0,$$
where $N'$ is the normal bundle of the isomorphic image of $\mathbb P^2\times\mathbb P^2$  under the 
projection of $\pi_P\colon\mathbb P^8\setminus\{P\}\to\mathbb P^7$ from 
a general point $P\in\mathbb P^8$. In particular, $N$ does not satisfy $\mathscr A_1$.\end{corollary}

\proof The first part follows from Theorem \ref{2n+1}, (2). To prove the second part we use the known fact that 
$X:=i(\mathbb P^2\times\mathbb P^2)$ is a Severi variety in $\mathbb P^8$, i.e. the projection of $\mathbb P^8$ from a general point of $\mathbb P^8$ maps $X$ biregularly onto a submanifolds $X'$ of $\mathbb P^7$. Then the conclusion
follows from Lemma \ref{r2}.\qed

\medskip

\begin{corollary}\label{ci2} Let $X$ be a submanifold of $\mathbb P^n$ of dimension $d\geq\frac{n}{2}$, with $n\geq 5$.  If $n=2d$ assume that $X$ is simply connected and $\mathscr O_X(1)$ is not divisible in $\Pic(X)$. If the normal bundle $N_{X|\mathbb P^n}$ is a direct sum of line bundles then $X$ is a complete intersection in 
$\mathbb P^n$. \end{corollary}

\proof Since $d\geq 3$ and $N_{X|\mathbb P^n}$ is a direct sum of line bundles,
Theorem \ref{2n+1}, (3) implies that $\Num(X)=\mathbb Z$. 
If $d>\frac{n}{2}$ then by Barth-Larsen theorem $X$ is simply connected and $\mathscr O_X(1)$ is not
divisible in $\Pic(X)$. If instead $n=2d$ we have this statements by hypotheses. It follows 
$\Num(X)=\Pic(X)=\mathbb Z[\mathscr O_X(1)]$.  Therefore in all cases the hypothesis that $N_{X|\mathbb P^n}$ is a direct sum of line bundles translates into $N_{X|\mathbb P^n}\cong\bigoplus\limits_{i=1}^{n-d}\mathscr O_X(a_i)$, with 
$a_i\geq 1$ and $n-d\leq\frac{n}{2}$. Then the conclusion follows from Theorem \ref{ci1} of Faltings.\qed

\begin{remarks}\label{russo} i) I am indebted to G. Ottaviani for calling my attention to the fact that if 
$d>\frac{n}{2}$, Corollary \ref{ci2} was first proved by Michael Schneider in \cite{[S]}. Although also based on Faltings' criterion of complete intersection (Theorem \ref{ci1} above), Schneider's proof is however different from ours because it uses the methods of \cite{[SGA2]} together with another result of Faltings \cite{[F1]} according to which every submanifold of $\mathbb P^n$ of dimension $>\frac{n}{2}$ satisfies the effective Grothendieck-Lefschetz condition 
$\Leff(\mathbb P^n,X)$.

\medskip

ii) Basili and Peskine proved that every nonsingular surface in $\mathbb P^4$ whose normal bundle is decomposable, is a complete intersection (see \cite{[BP]}). Their proof is based heavily on the methods developed in \cite{[EGPS]}. For a related result see M\'eguin \cite{[M]}. Partial results on the normal bundle of two-codimensional submanifolds of $\mathbb P^n$ (with $n\geq 5$) can be found in Ellia, Franco and Gruson 
\cite{[EF]}.\end{remarks}

\bibliographystyle{amsplain}

\end{document}